\documentclass[preprint,12pt]{elsarticle}

\usepackage{amssymb}
\usepackage{amsmath}
\journal{CNSNS}

\begin{document}

\begin{frontmatter}



\title{Neural network Approximations for Reaction-Diffusion Equations – Homogeneous Neumann Boundary Conditions and Long-time Integrations}


%\author[inst1]{Author One}

%\affiliation[inst1]{organization={Department One},%Department and Organization
%            addressline={Address One}, 
%            city={City One},
%            postcode={00000}, 
%            state={State One},
%            country={Country One}}

%\author[inst2]{Author Two}
%\author[inst1,inst2]{Author Three}

%\affiliation[inst2]{organization={Department Two},%Department and Organization
%            addressline={Address Two}, 
%            city={City Two},
 %           postcode={22222}, 
 %           state={State Two},
 %           country={Country Two}}

\author[1,2]{{Eddel El\'{i} Ojeda Avil\'{e}s}}
\ead{eddelojeda@gmail.com}

\author[2]{{Jae-Hun Jung}}
\ead{jung153@postech.ac.kr}
\author[1,2]{{Daniel Olmos Liceaga *}}\ead{daniel.olmos@unison.mx}
%\cormark[1]

\affiliation[1]{organization={Universidad de Sonora},
    addressline={Blvd. Luis Encinas y Rosales S/N}, 
    city={Hermosillo}, 
    state={Sonora},
    postcode={83000},
    country={Mexico}}

\affiliation[2]{organization={Pohang University of Science and Technology},
    %addressline={} 
    city={Pohang},
    % citysep={}, % Uncomment if no comma needed between city and postcode
    postcode={37673}, 
    %state={Korea},
    country={Korea}}

\cortext[cor1]{Corresponding author}

\begin{abstract}
Reaction-Diffusion systems arise in diverse areas of science and engineering. Due to the peculiar characteristics of such equations, analytic solutions are usually not available and numerical methods are the main tools for approximating the solutions. In the last decade, artificial neural networks have become an active area of development for solving partial differential equations. However, several challenges remain unresolved with these methods when applied to reaction-diffusion equations. In this work, we focus on two main problems. The implementation of homogeneous Neumann boundary conditions and long-time integrations. For the homogeneous Neumann boundary conditions, we explore four different neural network methods based on the PINN approach. For the long time integration in Reaction-Diffusion systems, we propose a domain splitting method in time and provide detailed comparisons between different implementations of no-flux boundary conditions. We show that the domain splitting method is crucial in the neural network approach, for long time integration in Reaction-Diffusion systems. We demonstrate numerically that domain splitting is essential for avoiding local minima, and the use of different boundary conditions further enhances the splitting technique by improving numerical approximations. To validate the proposed methods, we provide numerical examples for the Diffusion, the Bistable and the Barkley equations and provide a detailed discussion and comparisons of the proposed methods.

\end{abstract}




\end{frontmatter}


\section{Introduction}

Reaction-Diffusion systems are common in various research areas, such as physiology, ecology, chemistry, etc. \cite{brauer,epstein,keener_book,mu01,ty94}. Due to the non-linearities involved in the reaction, the common techniques for solving such equations are based on numerical methods, including finite differences, finite element, and pseudospectral methods \cite{Anguelov2005,Hu2012,os09,oj22,zhang2018}. While traditional methods continue to be developed and utilized, recent studies are increasingly shifting their focus towards machine learning approaches. The use of such techniques helps to better solve equations in high dimensions where traditional methods easily fail \cite{Beck2021_2,Blechschmidt2021,Han2018}; to predict solutions by training networks \cite{Li2020,ren2022,scholz2022}, and to speed up computations of traditional methods \cite{Mishra2018}. One of such methods uses machine learning techniques known as Physics informed Neural Network (PINN) \cite{raisi2018}. With this approach, instead of discretizing the equations, the problem is now to find a global minimum of a function named loss function which considers the given PDE, the initial, and boundary conditions. The use of such a technique is not completely new and arose about 30 years ago with the work of Lagaris \cite{lagaris98}, in which the trial solutions to solve ordinary differential equations were proposed such that initial conditions are satisfied exactly. The PINN approach has been widely used in different ways for RD equations. For example, \cite{Cho2021} used a specific form of a neural network to study the evolution of pulses with a constant traveling speed. In \cite{Mattey2022}, the authors developed a splitting in the time interval to deal with equations with high nonlinearities and high-order derivative operators. In \cite{giampaolo2022} the authors solved the Gray-Scott equations, a particular class of RD equations, with a modification of PINN where known data obtained with traditional numerical methods was used to obtain good solutions.

In this work we are interested in two different problems: (i) The solution with long time integrations
and (ii) the implementation of homogeneous Neumann boundary conditions. In order to deal with long-time integrations, we follow the idea of splitting the time interval as in Beck et. al. \cite{Beck2021_2,Beck2021}, Meng et al \cite{Meng2020} and Mattey and Ghosh \cite{Mattey2022}.  In \cite{Beck2021_2} the authors proposed a deep splitting method to solve parabolic PDEs in high dimensions by approximating the nonlinear part of the equation and changing it by a linear approximation. In Meng et al \cite{Meng2020}, the authors developed a method by implementing parallel computations, but in their work, it is not clear the efficiency of their method when the time integration is long. In \cite{Mattey2022}, the authors develop a method named bcPINN, in which the time interval is divided into smaller subintervals. In their work, they minimize a loss function that takes into account the current subinterval and all the previous optimized intervals. In our work, we take a completely different motivation from \cite{Beck2021_2} as we will consider the nonlinear problem for our computations. Also, in our work, we consider optimization over subdomains as in \cite{Mattey2022}, but in this case, we optimize in each subdomain sequentially. For each subdomain, we use a secondary network to pass information from the previous subdomain to the current one. \\

On the other hand, different authors have also worked on the implementation of boundary conditions. Reaction-diffusion equations are solved with Dirichlet boundary conditions \cite{Li2020,Zakeri2019}, periodic boundary conditions \cite{scholz2022} and/or mixed boundary conditions \cite{shekari2009}. Within the context of PINN, the numerical methods for reaction-diffusion equations have been mostly focused on the implementation of Dirichlet boundary conditions \cite{Cho2021,giampaolo2022} and periodic boundary conditions \cite{Mattey2022,raisi2018}. The implementation of the Neumann-type boundary conditions seems straightforward. However, it is not clear that the issue has been completely resolved, as mentioned in \cite{Cho2021}. The issue of satisfying Neumann or Dirichlet boundary conditions with the PINN approach may not be the problem. The present work focuses on solving reaction-diffusion equations with Neumann homogeneous boundary conditions. In this work, we propose three new different ways of implementing homogeneous Neumann boundary conditions: (i) Approximation of boundary conditions by finite difference methods, (ii) the Phase Field method and, (iii) the Mirror method. These developed techniques are compared with respect to the classical PINN.

In order to test our proposed methods, three equations are considered, i.e., the simple diffusion equation, the bistable equation and the Barkley's equations. Then, this paper is organized as follows. In Section \ref{Sec:PDE} we provide the basic introduction to the considered system of equations and the neural network methods. Then, in Section \ref{Sec:zfbc} we describe the methods we will test for implementing homogeneous Dirichlet boundary conditions. In Section \ref{Sec:subdomains} we discuss the subdomain approach used for long-time integration dynamics. After this, we present our results in Section \ref{Sec:Results}. Finally, we conclude with a section of conclusions and discussions (Sec. \ref{Sec:Conc_Disc}). 

\section{Solution of partial differential equations}\label{Sec:PDE}
Consider the following reaction-diffusion set equations
\begin{equation}\label{Eq:RD_PDERn}
    \frac{\partial U^i}{\partial t}=\nabla \cdot \mathbf{D}^i(X)\nabla U^i +F(U,t),
\end{equation}
for $i=1,...,d_{out}$, $U=(U^1,U^2,\cdots U^{d_{out}})$, $X \in \mathbf{R}^n$ and such that the solution is sought over a connected domain $\Omega=\Omega_x \times \Omega_t \subset \mathbf{R}^n \times \mathbf{R}$, with initial condition $U(x,0)=U_0(x)$ and boundary conditions $\nabla U^i\cdot \eta=0$ at the boundary of $\Omega_x$ with $\eta$ being the perpendicular vector to the boundary. In Eq. \ref{Eq:RD_PDERn}, $U^i: \mathbf{R}^n \times \mathbf{R} \rightarrow \mathbf{R}^n$, $\nabla$ represents the gradient operator, $\mathbf{D}^i(X)$ is the medium conductivity matrix for the $i$th species, and $F: \mathcal{R}^n \times \mathcal{R} \rightarrow \mathcal{R}^n$ is the reaction term. When $\mathbf{D}^i(X)=\mathbf{D}^i$ is a constant matrix, the equation becomes
\begin{equation}
    \frac{\partial U^i}{\partial t}=\mathbf{D}^i\nabla^2 U^i +F(U,t),
\end{equation}
for $i=1,...,d_{out}$.\\
\begin{figure}[h]
	\centering
	\includegraphics[width=12cm, height=5.5cm]{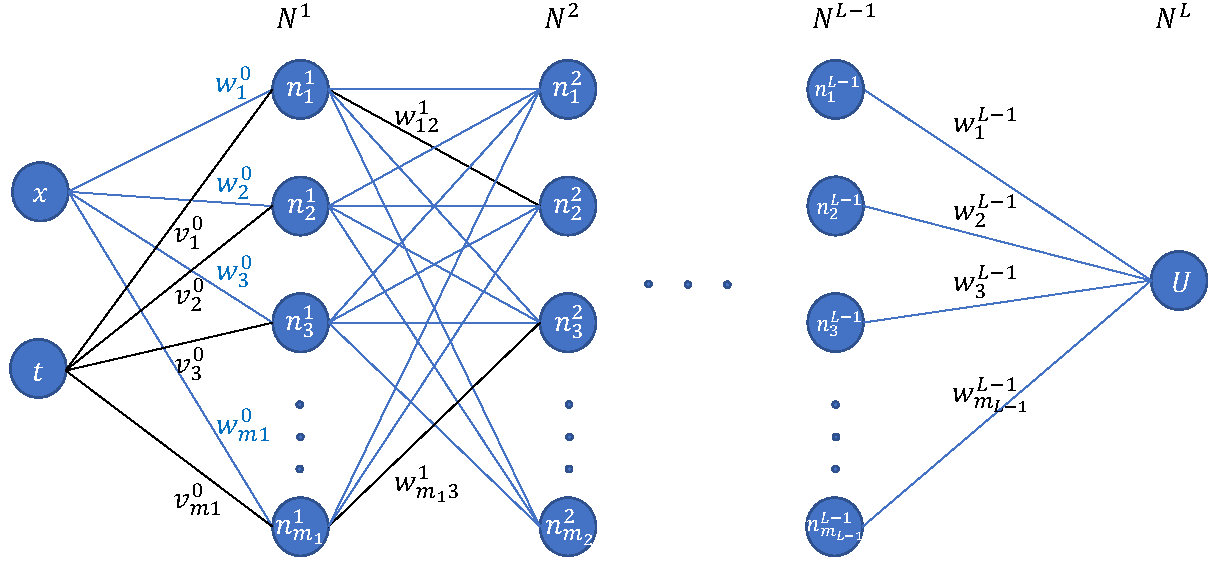}
	\caption{A neural network of $L$ layers with an input of $x$ and $t$.}
	\label{Fig:NN_formulation}
\end{figure}
In this work we focus on one-dimensional problems, i.e. 
\begin{equation}
    \frac{\partial U^i}{\partial t}=D^i\frac{\partial^2 U^i}{\partial x^2} +F(U,t),
\end{equation}
for $i=1,...,d_{out}$, with the initial and boundary conditions, 
where $D^i$ is the diffusion coefficient for the $i$th species, and the solution is sought over the domain $\Omega=\Omega_x\times\Omega_t=\{(x,t)|x\in\Omega_x=[a,b], a<b, t\in\Omega_t=[t_0,T]\}$. In this work we have a particular interest in no-flux boundary conditions, i.e.
\begin{equation}
    \left.\frac{\partial U}{\partial x}\right\vert_{x=a}=0, \qquad\qquad\qquad\qquad
    \left.\frac{\partial U}{\partial x}\right\vert_{x=b}=0.
\end{equation}
%%%%%%%%%%%%%%%%%%%%%%%%%%%%%%%%%%%%%%%%%%%%%%%%%%%%%%%%%%%%%%%%%%%%
We seek a trial solution $U_T$ at $(x,t)$ as shown in Fig. \ref{Fig:NN_formulation} and. It is given by
\begin{equation}
    U_T(x,t)=\sum_{i=1}^{m_{L-1}}w_i^{L-1}n_i^{L-1}+b^{L-1}
\end{equation}
where
$$
n_i^{L-1}=\sigma\left(\sum_{r=1}^{m_{L-2}}w_{ri}^{L-2}n_r^{L-2}+b_i^{L-2}\right)
$$
where $\sigma$ is the activation function applied to each vector component and
\begin{equation}
    n_i^j=\sigma\left(\sum_{r=1}^{m_{j-1}}w_{ri}^{j-1}n_r^{j-1}+b_i^{j-1}\right)
\end{equation}
for $j=2,...,L-1$. In the case of the values $n_i^1$, $i=1,...,m_1$, we simply have
$$n_i^1=\sigma(w_i^0x+v_i^0t+b_i^0).$$ 
This formulation can be stated in a compact form as given in \cite{Lu2021} as follows
$$
\begin{array}{rl}
\textrm{input layer:}& N^0=(x,t),\\
\textrm{hidden layers:} & N^l(x,t)=\sigma(W^{l-1} N^{l-1}(x,t)+b^{l-1}) \in \mathbf{R}^{N_l} \qquad \textrm{for} \qquad 1\leq l\leq L-1,\\
\textrm{output layer:}& U_T(x,t)=N^L(x,t)=W^{L-1}N^{L-1}(x,t)+b^{L-1} \in \mathbf{R}^{d_{out}}
\end{array}
$$
where $\mathbf{R}^{d_{out}}$ is the number of dependent variables of the system.
With the trial function given above, we minimize the loss function defined by
\begin{equation}
\begin{aligned}
    L(p)=&\frac{1}{N_{int}}\sum_{j=1}^{N_{int}} \sum_{i=1}^{d_{out}} \left.\left( (U^i_T(p,x_j,t_j))_t - D^i(U^i_T(p,x_j,t_j))_{xx} - F^i(U_T(p,x_j,t_j),t_j) \right)^2\right\vert_{(x_j,t_j)\in \Omega^o}\\
    & +\frac{1}{N_{bc}}\sum_{j=1}^{N_{bc}}\sum_{i=1}^{d_{out}}
    \left.\left(
    (U^i_T(p,x_j,t_j))_x \right)^2\right\vert_{(x_j,t_j)\in \partial \Omega_x \times \Omega_t}\\
     &+\frac{1}{N_{ic}}\sum_{j=1}^{N_{ic}}\sum_{i=1}^{d_{out}}\left.\left( U^i_T(p,x_j,0)-U^i_0(x_j) \right)^2\right\vert_{x_j\in\Omega_x}, 
\end{aligned}    
\end{equation}
with respect to the weights and biases, denoted by $p$. In the equation, $\Omega^o$ and  $\partial \Omega$ refer to the interior and the boundary of the $\Omega$ set, respectively. The superscript $i$ refers to the $i$th component of $U_T$ with $d_{out}$ components, and $\lVert\cdot\rVert$ denotes the vector $L_2$ norm; $N_{int}$, $N_{bc}$ and $N_{ic}$ are the number of random sampling points in the interior, the spatial boundary and at $t=0$, respectively. The subscripts $t,x$ and $xx$ denote the partial derivatives of function $U^i_T$ obtained by automatic differentiation.

\section{Implementing homogeneous Neumann boundary conditions}\label{Sec:zfbc}

One of the main concerns in this paper is an appropriate implementation of no-flux boundary conditions. Implementing homogeneous Neumann boundary conditions is straightforward based on the PINN approach. However, as stated in \cite{Cho2021}, solving PDEs with Neumann boundary conditions for Reaction-Diffusion equations by NN, is not a completely resolved issue as it seems to provide errors in the solution. Therefore, in the present section, we describe three different ways to implement homogeneous Neumann boundary conditions additionally to PINN in order to improve the calculation of the PDE solution. These strategies are:  1) Use of Finite Difference methods to approximate the derivative at the boundary 2) The Phase Field method and 3) The Mirror trick. 

\subsection{Finite difference boundary condition}

One of the issues of implementing boundary conditions, is that such conditions are imposed only at a single point. Therefore, one idea is to increase the number of points used to impose the boundary conditions. This can be done by increasing the number of sampling points at the boundary or approximate the boundary condition such that
points in the interior of the domain are also involved in the minimization. We can apply no-flux boundary conditions by using a finite difference approximation at the boundary. 

\subsubsection{Second order approximation.} 

Consider the left and right boundaries and approximate the derivative with its second-order approximation given by
\begin{equation}\label{Eq:der_2ordl}
    f'(x)\approx \frac{-3f(x)+4f(x+\Delta x)-f(x+2\Delta x)}{2\Delta x},
\end{equation}
and
\begin{equation}\label{Eq:der_2ordr}
    f'(x)\approx \frac{3f(x)-4f(x-\Delta x)+f(x-2\Delta x)}{2\Delta x},
\end{equation}
for the left and right boundaries. With this in mind, we then minimize the loss function given by
\begin{equation}
\small 
\begin{aligned}
    L(p)=L_{int}(p)+L_{ic}(p)+L_{bc}(p),
\end{aligned}
\end{equation}
where $L_{int}$ is the loss function defined in the domain interior, $L_{ic}$ refers to the loss function due to the initial condition, and $L_{bc}$ refers to the loss at the boundary. These loss functions are given by
\begin{equation}
\small
\begin{aligned}
    L_{int}(p)=&\frac{1}{N_{int}}\sum_{j=1}^{N_{int}} \sum_{i=1}^{d_{out}} \left.\left( (U^i_T(x_j,t_j))_t - d_{ii}(U^i_T(x_j,t_j))_{xx} - F^i(U_T(x_j,t_j),t_j) \right)^2\right\vert_{(x_j,t_j)\in \Omega^o},\\
\end{aligned}    
\end{equation}

\begin{equation}
\small
\begin{aligned}
    L_{ic}(p)= & \frac{1}{N_{ic}}\sum_{j=1}^{N_{ic}}\sum_{i=1}^{d_{out}}\left.\left( U^i_T(x_j,0)-U^i_0(x_j) \right)^2\right\vert_{t=0},
\end{aligned}    
\end{equation}
and

\begin{equation}
   \small
\begin{aligned}
   & L_{bc}(p)=\\
    & \frac{1}{N_{bc}}\sum_{j=1}^{N_{bc}} \sum_{i=1}^{d_{out}} \left.\left( (U^i_T(x_L^1,t_j))_t - D^i(U^i_T(x_L^1,t_j))_{xx} - F^i(U_T(x_L^1,t_j),t_j) \right)^2\right\vert_{t\in \Omega_t}+\\
     & \frac{1}{N_{bc}}\sum_{j=1}^{N_{bc}} \sum_{i=1}^{d_{out}} \left.\left( (U^i_T(x_L^2,t_j))_t - D^i(U^i_T(x_L^2,t_j))_{xx} - F^i(U_T(x_L^2,t_j),t_j) \right)^2\right\vert_{t\in \Omega_t}+\\
         & \frac{1}{N_{bc}}\sum_{j=1}^{N_{bc}} \sum_{i=1}^{d_{out}} \left.\left( (U^i_T(x_R^{-1},t_j))_t - D^i(U^i_T(x_R^{-1},t_j))_{xx} - F^i(U_T(x_R^{-1},t_j),t_j) \right)^2\right\vert_{t\in \Omega_t}+\\
     & \frac{1}{N_{bc}}\sum_{j=1}^{N_{bc}} \sum_{i=1}^{d_{out}} \left.\left( (U^i_T(x_R^{-2},t_j))_t - D^i(U^i_T(x_R^{-2},t_j))_{xx} - F^i(U_T(x_R^{-2},t_j),t_j) \right)^2\right\vert_{t\in \Omega_t}+\\
          & \frac{1}{N_{bc}}\sum_{j=1}^{N_{bc}} \sum_{i=1}^{d_{out}} \left.\left( 3U^i_T(x_L,t_j)  -4U^i_T(x_L^1,t_j)+U^i_T(x_L^{2},t_j) \right)^2\right\vert_{t\in \Omega_t}+\\
        & \frac{1}{N_{bc}}\sum_{j=1}^{N_{bc}} \sum_{i=1}^{d_{out}} \left.\left( 3U^i_T(x_R,t_j)  -4U^i_T(x_R^{-1},t_j)+U^i_T(x_R^{-2},t_j) \right)^2\right\vert_{t\in \Omega_t}.\\
\end{aligned}    
\end{equation}
where $x_L^1,x_L^2,x_R^{-1},x_R^{-2}$ are equal to $x_L+\Delta x, x_L+2\Delta x, x_R-\Delta x$ and $x_R-2\Delta x$, respectively. In the formulation $U_T=(U_T^1,U_T^1,...,U_T^{d_{out}})$. The first four terms correspond to the interior points used to enforce the boundary conditions and must satisfy the governing partial differential equation (PDE). The last two terms are included to ensure compliance with Eqns. \ref{Eq:der_2ordl} and \ref{Eq:der_2ordr} at the boundary. We have that the equation is satisfied at the interior points and initial conditions. Boundary conditions are satisfied approximately using a finite difference (FD) formulation. The main advantage is that boundary conditions involve points in the interior. An important parameter is the value of $\Delta x$, which, by taking it small enough, will provide a better approximation of the derivative.\\

Similar formulations are obtained for higher-order approximations of the derivative. For example, by using the third-order derivative approximation of the first derivative of a function at a boundary we have
$$
    f'(x)\approx \frac{-11f(x)+18f(x+\Delta x)-9f(x+2\Delta x)+2f(x+3\Delta x)}{6\Delta x},
$$
and
$$
    f'(x)\approx \frac{11f(x)-18f(x-\Delta x)+9f(x-2\Delta x)-2f(x-3\Delta x)}{6\Delta x},
$$
for the right and left boundaries, respectively. 
Clearly, we can use higher-order approximations, up to having spectral accuracy by using Chebyshev polynomials. 

\subsection{The Phase Field method}

This method is used to solve partial differential equations of the reaction-diffusion type on irregular domains. In \cite{Bueno2006}, the authors developed this method and used a Fourier expansion for the general formulation of the RD equations given by

\begin{equation*}%\label{Eq:RD_PDERn}
    \frac{\partial U^i}{\partial t}=\nabla \cdot \mathbf{D}^i(X)\nabla U^i +F(U,t),
\end{equation*}

 with initial condition $U^i(x,0)=U_{i0}(x)$ and no flux boundary conditions
 $$
 \vec{n}\cdot\mathbf{D}^{i}\nabla U^i=0\qquad \textrm{on}\qquad \partial\Omega, 
 $$
where $\vec{n}$ is the normal vector to the boundary $\partial\Omega$.
In this case, $\mathbf{D}^{j}$, $j=1,...,d_{out}$ is a diffusion coefficient matrix that may depend on space. The idea of solving a RD equation on irregular domains, implies to solve the extended problem
 $$
\partial_t(\phi^{(\xi)}U^{i(\xi)})=\nabla\cdot(\phi^{(\xi)}\mathbf{D}^{i}\nabla U^{i(\xi)})+\phi^{(\xi)}f(U^{(\xi)},t),
 $$
 for $U^{i(\xi)}$ on an enlarged domain $\Omega'$, such that
 (i) $\Omega\subset \Omega'$ and (ii) $\partial\Omega \cap \partial \Omega'=\emptyset$. $\phi^{(\xi)}$ is continuous in $\Omega'$, takes the value of 1 inside $\Omega$, smoothly decays to zero outside $\Omega$, with $\xi$ identifying the width of the decay (See \cite{Bueno2006} for details). 
One advantage of the Phase Field method with the NN technique is that it is possible to apply the method in two ways 
\begin{enumerate}
    \item Take a continuous approximation of $X_{\Omega}$. In this case, we have to minimize

$$
\small
\begin{aligned}
    L(p)=&\frac{1}{N_{int}}\sum_{j=1}^{N_{int}} \sum_{i=1}^{d_{out}} \left.\left(\phi^{(\xi)} (U^i_T(x_j,t_j))_t - d_{ii}(\phi^{(\xi)}(U^i_T)_x(x_j,t_j))_{x} - \phi^{(\xi)}F^i(U_T(x_j,t_j),t_j) \right)^2\right\vert_{(x_j,t_j)\in \Omega'^{o}}+\\
    & \frac{1}{N_{bc}}\sum_{j=1}^{N_{bc}}\sum_{i=1}^{d_{out}}
    \left.\left(
    BC(U^i_T(x_j,t_j)) \right)^2\right\vert_{(x_j,t_j)\in \partial \Omega_x\times\Omega_t}
     +\frac{1}{N_{ic}}\sum_{j=1}^{N_{ic}}\sum_{i=1}^{d_{out}}\left.\left( U^i_T(x_j,0)-U^i_0(x_j) \right)^2\right\vert_{x_j\in\Omega_x}. 
\end{aligned}    
$$
   
    \item Take $X_{\Omega}$ directly. Here we need to minimize
$$
\small
\begin{aligned}
    L(p)=&\frac{1}{N_{int}}\sum_{j=1}^{N_{int}} \sum_{i=1}^{d_{out}} \left.\left( (U^i_T(x_j,t_j))_t - d_{ii}(U^i_T(x_j,t_j))_{xx} - F^i(U_T(x_j,t_j),t_j) \right)^2\right\vert_{(x_j,t_j)\in \Omega^o}+\\
    & \frac{1}{N_{bc}}\sum_{j=1}^{N_{bc}}\sum_{i=1}^{d_{out}}
    \left.\left(
    (U^i_T(x_j,t_j))_x \right)^2\right\vert_{(x,t)\in \partial \Omega_x'}
     +\frac{1}{N_{ic}}\sum_{j=1}^{N_{ic}}\sum_{i=1}^{d_{out}}\left.\left( U^i_T(x_j,0)-U^i_0(x_j) \right)^2\right\vert_{x_j\in\Omega_x}. 
\end{aligned} 
$$
Observe that in this case, the region $\Omega'\cap\Omega^c$ does not play any role in the minimization process.
\end{enumerate}
The term involving the boundary conditions in the first formulation may consider Dirichlet boundary conditions or periodic boundary conditions. In the case of Dirichlet boundary conditions the values to be taken at the boundary, are given by the appropriate asymptotically stable steady states of the system.

\subsection{The Mirror technique}
\begin{figure}[h]
	\centering
	\includegraphics[width=11cm, height=7.0cm]{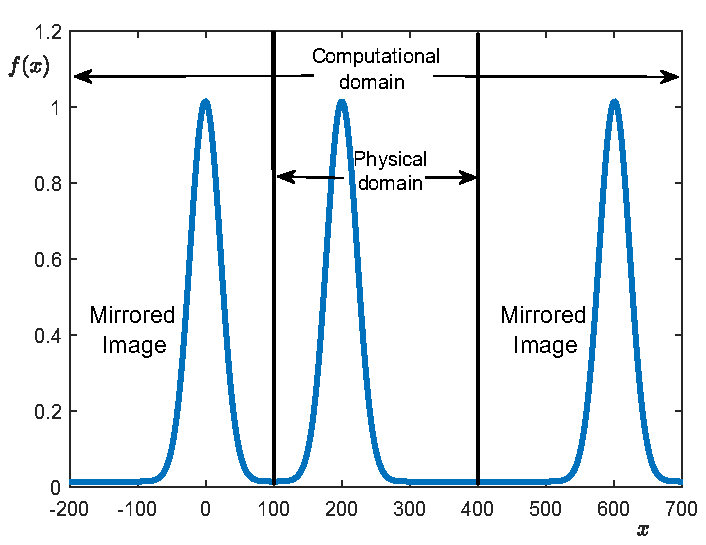}
	\caption{Physical and Computational domains for the Mirror technique. The mirrors are located at $x=100$ and $x=400$.}
	\label{Fig:mirror_method}
\end{figure}

This method is based on the physical properties of diffusion. In order to implement no-flux boundary conditions, we extend as in the Phase Field method, our computational domain. For our dimensional physical domain $\Omega=\Omega_x\times\Omega_t =\{(x,t)|x\in\Omega_x=[x_L,x_R], x_L<x_R, t\in\Omega_t=[t_0,T]\}$, we consider the extended domain  $\Omega'=\Omega'_x\times\Omega_t=\{(x,t)|x\in\Omega'_x=[x_L,x_R], x_L<x_R, t\in\Omega_t=[t_0,T]\}$. For a function $f(x)$ that is defined on $\Omega_x$, we define its Mirror extension $\hat{f}$ over the domain $\Omega'_x$ as
$$
\hat{f}(x)=
\left\{
\begin{array}{lll}
f(2x_L-x)& \textrm{for} & x \in [x_L-L,x_L]\\
f(x)&\textrm{for} & x \in [x_L,x_R]\\
f(2x_R-x)& \textrm{for} & x \in [x_R,x_R+L].\\
\end{array}
\right.
$$
From Fig. \ref{Fig:mirror_method}, we can see that the reflected images at $x=x_L$ and $x_R$ create the same effect as no flux boundary conditions at the physical boundaries. This is due to the fact that the front has an exact copy on the other side which do not allows the diffusion towards each other side as there is no gradient at the physical boundaries ($x=100$ and $x=400$). 

The procedure of implementing the extended mirrored domain, for the NN procedure begins by extending the initial condition over the computational domain. Therefore, minimization of the problem takes place over the larger domain with Dirichlet boundary conditions, i.e. we look to minimize

\begin{equation}\label{Eq:lossMIRROR}
\begin{aligned}
    L(p)=&\frac{1}{N_{int}}\sum_{j=1}^{N_{int}} \sum_{i=1}^{d_{out}} \left.\left( (U^i_T(p,x_j,t_j))_t - D^i(U^i_T(p,x_j,t_j))_{xx} - F^i(U_T(p,x_j,t_j),t_j) \right)^2\right\vert_{(x_j,t_j)\in \Omega'^o}+\\
    & \frac{1}{N_{bc}}\sum_{j=1}^{N_{bc}}\sum_{i=1}^{d_{out}}
    \left.\left(
    U^i_T(p,x_j,t_j)-U_{\infty}^i \right)^2\right\vert_{(x_j,t_j)\in \partial \Omega'_x \times \Omega_t}\\
     &+\frac{1}{N_{ic}}\sum_{j=1}^{N_{ic}}\sum_{i=1}^{d_{out}}\left.\left( U^i_T(p,x_j,0)-U^i_0(x_j) \right)^2\right\vert_{x_j\in\Omega'_x}, 
\end{aligned}    
\end{equation}
where at the boundary of the larger domain, Dirichlet boundary conditions are imposed and the value to be taken is $U^i(x_L,t)=U^i(x_R,t)=U^i_{\infty}$, which corresponds to the coordinates of the steady-state values obtained from the system without diffusion. 

An important issue arises due to Dirichlet boundary conditions at the extended computational domain. If we make a mirror image only when $t=0$, the Dirichlet boundary conditions over the extended computational domain will create an unbalance in the equations, giving in the long integration time a break up in the solution's symmetry and therefore a break up in the no-flux boundary conditions at the physical domain. In order to overcome this problem, a natural way to avoid this problem is to separate the time domain into subdomains, a technique to be discussed in the next section. In this case, every time we change the domain, we redefine the initial condition over the extended domain.\\

\section{Separating in subdomains: Long time integrations}\label{Sec:subdomains}
 One important task when solving the reaction-diffusion equations is to consider long time integration. In these equations, simulations over long periods are needed, as there are transients that may take long to fade away. However, for a given spatial domain the larger the integration time, the larger the amount of sampling points over the domain is needed to preserve point density. For large integration times, the required number of points to maintain the sampling density may increase, leading to memory problems. Also, if we focus only on ordinary differential equations and a single hidden layer to approximate a periodic solution, it is clear that as we increase the integration time, it will need to use more nodes to approximate appropriately the solution. 
 
 A way to overcome this problem is by dividing the complete domain $\Omega=\{(x,t)/x\in\Omega_x=[a,b],t\in\Omega_t=[0,T]\}$, into smaller subdomains. The procedure is as follows

\begin{enumerate}
\item  Consider the problem to solve the RD equation over the domain $\Omega=\{(x,t)/x\in\Omega_x,t\in\Omega_t\}$ with initial and boundary conditions. 
\item Consider the domains sequence given by $\Omega_k=\{(x,t)/x\in\Omega_x,t\in \Omega_t^k=[t_{k-1},t_k]\}$, $k=1,..,N$, where $t_0<t_1<t_2<\cdots<t_N=T$. Also, consider two neural networks, $NN_{IC}^k$ and $NN_{PDE}^k$ for each subdomain. $NN_{IC}^k$ is used to approximate the initial condition of the equation at the $\Omega_k$ subdomain and $NN_{PDE}^k$ is used to approximate the solution of the Dirichlet or Neumann problem on the $\Omega_k$ domain.
\item For the first subdomain, approximate the initial condition with $NN_{IC}^1$ and obtain the optimized parameters denoted by $w_{IC}^1$. This approximation is denoted by $U^{w_{IC}^1}$.
\item For $\Omega_1=\{(x,t)/x\in\Omega_x,t\in\Omega_t^1\}$, the loss function to be minimized is

\begin{equation}\label{Eq.Losstotal}
\small
\begin{aligned}
Loss_{Total_1}(p)=&Loss_{PDE_1}+Loss_{BC_1}+Loss_{IC_1}\\
    =&\frac{1}{N_{int}^1}\sum_{j=1}^{N_{int}^1} \sum_{i=1}^{d_{out}} \left.\left( (U^i_T(x_j,t_j))_t - d_{ii}(U^i_T(x_j,t_j))_{xx} - F^i(U_T(x_j,t_j),t_j) \right)^2\right\vert_{(x,t)\in \Omega_1^o}+\\
    & \frac{1}{N_{bc}^1}\sum_{j=1}^{N_{bc}^1}\sum_{i=1}^{d_{out}}
    \left.\left(
    (U^i_T(x_j,t_j))_x \right)^2\right\vert_{x\in \partial \Omega_x, t \in \Omega_t^1}
     +\frac{1}{N_{ic}^1}\sum_{j=1}^{N_{ic}^1}\sum_{i=1}^{d_{out}}\left.\left( U^i_T(x_j,0)-U^{w_{IC}^1} \right)^2\right\vert_{t=\min\Omega_t^1}. 
\end{aligned}    
\end{equation}
%%%%%%%%%%%%%%%%%%%%%%%%%%%%%%%%%%%%%%%%%%%%%%%
\begin{figure}[h]
	\centering
\includegraphics[width=14cm, height=7.0cm]{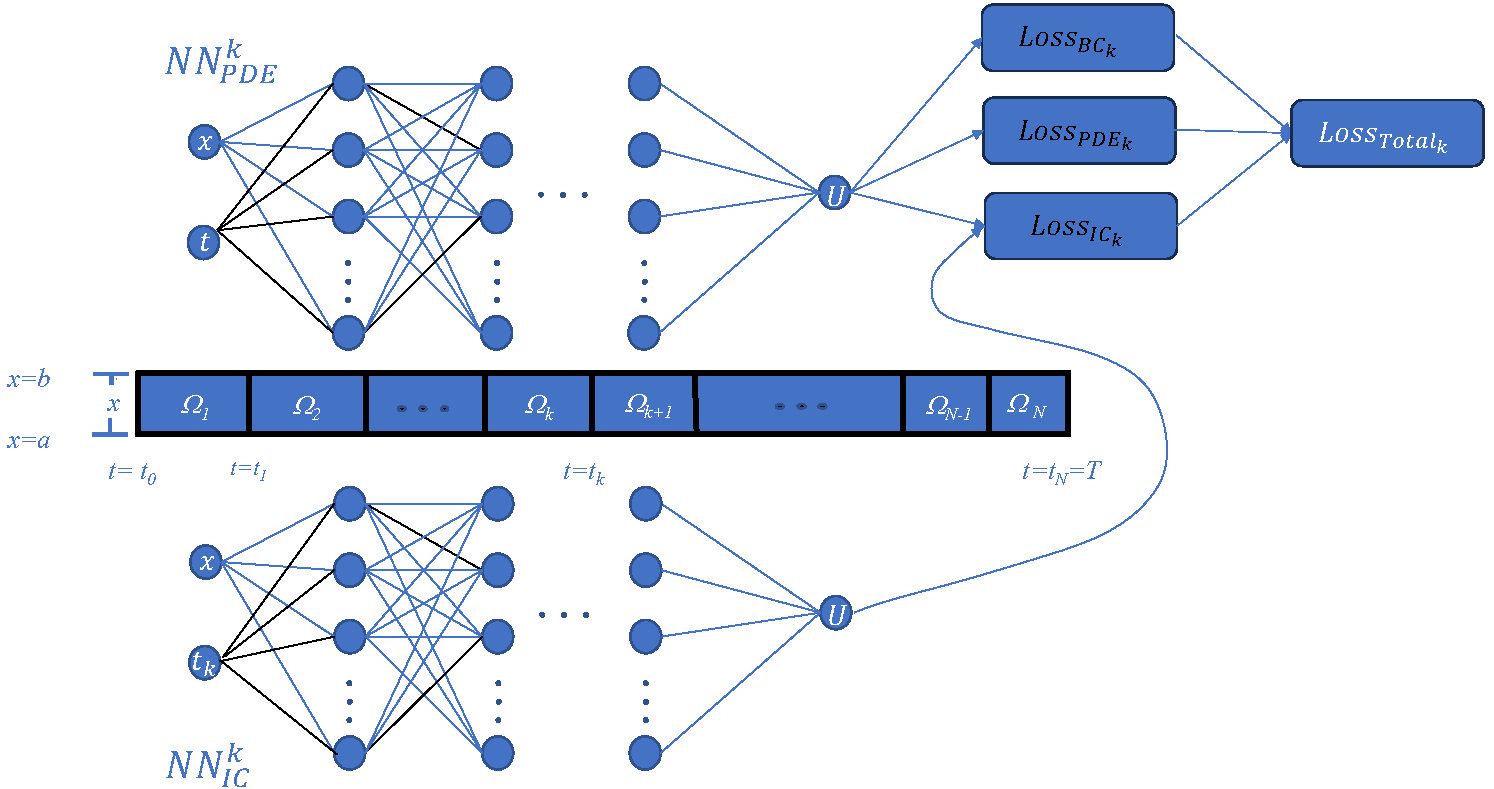}
	\caption{Configuration of the two neural networks involved in the minimization of $Loss_{Total}$ (Eq. \ref{Eq.Losstotal}) for each subdomain.}
	\label{Fig:NN_2N}
\end{figure}
%%%%%%%%%%%%%%%%%%%%%%%%%%%%%%%%%%%%%%%%%%%%%%%
In this case, $t_1\leq T$. Observe that the last term of this sum corresponds to approximate a known function, given by $NN_{IC}^1$ evaluated at $t=t_0$. Eq. \ref{Eq.Losstotal} is optimized when the average loss in $p$ iterations is less than an $\epsilon$ value given. A visualization on how the loss function is calculated is shown in Fig. \ref{Fig:NN_2N}.

\item Once the solution at $\Omega_k$ is obtained, we use the trained weights $w_{INT}^k$ to basically evaluate any vector of the form $U(x,t_k)$, i.e. an approximated solution at $t=t_k$ is used as a new initial condition for the time interval $[t_k, t_{k+1}]$. In this case, we take $w_{IC}^{k+1}$ as the optimized weights at $\Omega_k$, i.e. equal to $w_{INT}^k$.  

\item The procedure repeats until, the equation is solved for $t=T$.
\end{enumerate}

This method, takes the advantage of the optimized weights at domain $\Omega_k$ to optimize the weights at domain $\Omega_{k+1}$. This is a great advantage when optimizing the weights at $\Omega_{k+1}$ as the initial weights search contains information from the previous domain and the initial search may converge faster than the usage of a random or predetermined weights.

The use of this subdomain technique seems to provide advantages. Aside from the memory saving, one can think on the approximation of a function $f$ with a given $U_T$, over a domain $\Omega$. Clearly, if $||U_T-f||=\delta <<1$ for $(x,t)\in \Omega$, then  $||U_T-f||\leq\delta$ for $(x,t) \in \hat{\Omega} \subset \Omega$. Therefore, one would expect a better result when the minimization domain is reduced.

%\subsection{Variable learning rate}
In order to solve our equations, the algorithm uses the Adam method with variable learning rate. For this, we define the loss thresholds $\epsilon_1$, $\epsilon_2$,..., $\epsilon_n$ and learning rates $\eta_0$, $\eta_1$,..., $\eta_n$. Initially, Adam method uses a learning rate of $\eta_0$. As soon as the value of the loss function is less than $\epsilon_k$, the learning rate takes the value of $\eta_k$, even that consequent values for the loss function are larger than $\epsilon_k$.

%\subsection{Changing subdomain criteria}
One of the most important issues in this method is to provide criteria to change the optimization in subdomain $\Omega_j$ to subdomain $\Omega_{j+1}$. This criteria has to assure that the solution in $\Omega_j$ is a reliable one, before moving to the next subdomain. Clearly, these criteria must work as a stopping condition when only one domain is considered.

In this work we consider three different criteria: (i) The first is to consider a function $A_L$ as the average of the last $N$ values of the loss function. When $A_L$ becomes less than a predetermined value $\epsilon_L$ it can be considered adequate to change the domain. However, this criteria will not necessarily work for multiple subdomains. As it will be shown in the examples, for each subdomain, the $\epsilon_L$ value required to reach an acceptable converged solution, may change from interval to interval. (ii) The second is the evaluation of the loss function at epochs $s$ and $s+p$. Then, we calculate the quantity
\begin{equation}\label{Eq_Dl}
D_L=|\log_{10}(Loss(s+p))-\log_{10}(Loss(s))|.
\end{equation}
when this quantity becomes less than a given value $\epsilon_L$, then we assume the loss function will not change considerably after $p$ epochs and therefore, a change to the next subdomain is taken. (iii) The third criterion takes into account that the loss function can have a large amount of fluctuations. In this situation, we define the vector $L_M$ where $L_M[s]$ has the minimum value of the loss function up to epoch $s$. From here, we can find the change this function has after $p$ epochs. By using the logarithmic scale to measure this change, we look to satisfy the condition $S_L<\epsilon_L$, where
\begin{equation}\label{Eq_Sl}
S_L=|\log_{10}(L_M[s])-\log_{10}(L_M[s-p])|.
\end{equation}
When there is not a new minimum value found after $p$ epochs, we conclude the method will not be able to reduce more the loss function and a change of domain follows.

The $A_L$ value is used to approximate functions or when single domain computations are searched, $D_L$ is considered when changing subdomains and the loss function does not have significant fluctuations and $S_L$ can be used in any situation. 

\section{Numerical Results}\label{Sec:Results}
 In order to see how our method performs, we consider three equations: i) the diffusion equation, ii) the bistable equation and, iii) the Barkley's equations. 
 
\subsection{Diffusion Equation}
  
The first equation to be solved is the diffusion equation, which is given by 
\begin{equation}
    u_t=Du_{xx},
\end{equation}
for $x\in[-3,3]$, $t\in[0,3]$, $D=1$, initial condition $u(x,0)=u_0(x)=1+\cos\left(\frac{\pi}{3}x\right)$ and no-flux boundary conditions, i.e. $\frac{\partial u}{\partial x}|_{x=x_L}=\frac{\partial u}{\partial x}|_{x=x_R}=0$. This equation can be solved exactly and its solution is
$$
u(x,t)=1+\cos\left(\frac{\pi x}{3}\right)e^{-\left(\frac{\pi}{3}\right)^2t}.
$$
\begin{figure}[hbt!]
	\centering
	\includegraphics[width=13cm, height=15.0cm]{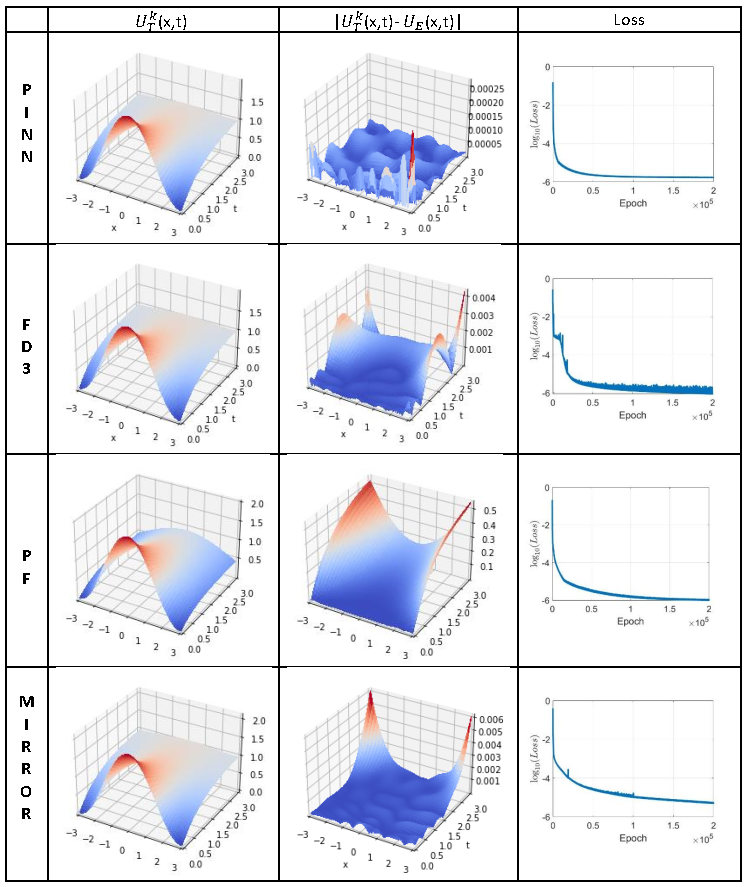}
	\caption{Performance of the four methods to solve the diffusion equation with no-flux boundary conditions. The first column, shows the approximation of the solution obtained with NN, for a particular epoch $k$. The second column, shows the error at each point of the mesh when compared with the exact solution. This is to show where the main source of errors are located. The last column shows the evolution of the loss function respect to the epoch.}
	\label{Fig:dif_icbcapprox}
\end{figure}

\begin{figure}[h]
	\centering
	\includegraphics[width=13cm, height=5.5cm]{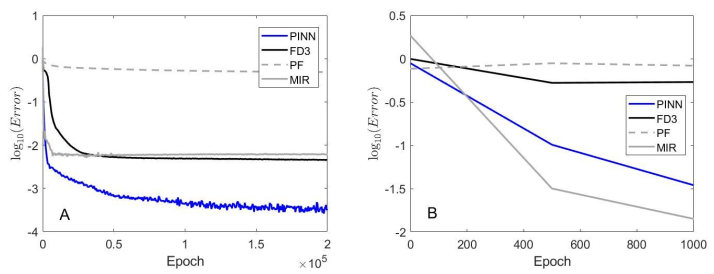}
	\caption{Obtained error for the different methods to solve the diffusion equation with no-flux boundary conditions. (A) the  $L_{\infty}$ error versus the epoch, (B) A zoom for the first 1000 epochs. The error is obtained over a regular mesh of 400 points in the $x$ direction and 250 points in the $t$ direction.} 
	\label{Fig:error_difussion}
\end{figure}

The results of implementing the methods are summarized on Fig. \ref{Fig:dif_icbcapprox} for $2\times10(5)$ epochs. In the figure, we show the approximation of the solution, the point-wise error in absolute value respect to the exact solution, and the loss as a function of the epochs. The results in the figure, clearly show that PINN, FD3 and Mirror provide with a better solution compared to the Phase Field method. The second and third columns show a clear advantage of the PINN and the Mirror methods over the other two. The Mirror method suggest an improvement far from $t=3$ compared to the PINN method. However if an $l_{\infty}$ measurement is taken, the PINN method performs better. From the third column it appears that also, PINN method converges faster to the exact solution than the Mirror method. This is corroborated on Fig. \ref{Fig:error_difussion}, where the $l_{\infty}$ error over a given mesh is computed. Clearly, the PINN method converges faster. However, it is worth to mention that in order to obtain an almost indistinguishable approximation to the exact solution, for the PINN method, it is required to run 1000 epoch, whereas for the Mirror method it only suffices to apply 500 epoch (Fig \ref{Fig:error_difussion}B).

On the other hand, in an attempt to improve the performance of the methods based on finite difference approximations, additionally to the FD3 method, the FD4 and FD2 methods were implemented as well. However, best results were acquired with FD3. The value of $\Delta_x=0.05$ was used for all the finite difference methods as it provided better results compared to other values. Also, for FD3 method, the separation in subdomains was considered, giving a reduction of about one order of magnitude. However, the PINN method in a single domain, performed better than this improvement of FD3. 

\begin{table}[hbt!]
\center
\begin{tabular}{|l|l|l|l|l|}
\hline
                & \multicolumn{1}{c|}{PINN}                                   & \multicolumn{1}{c|}{FD3}                                    & \multicolumn{1}{c|}{PF}                                     & \multicolumn{1}{c|}{Mirror}                                 \\ \hline
Network         & \begin{tabular}[c]{@{}l@{}}4 Layers\\ 60 Nodes\end{tabular} & \begin{tabular}[c]{@{}l@{}}4 Layers\\ 60 Nodes\end{tabular} & \begin{tabular}[c]{@{}l@{}}4 Layers\\ 60 Nodes\end{tabular} & \begin{tabular}[c]{@{}l@{}}4 Layers\\ 60 Nodes\end{tabular} \\ \hline
$N_{int}$       & 8000                                                        & 8000                                                        & 8000                                                        & 24000                                                       \\ \hline
$N_{bc}$        & 10000                                                       & 10000                                                       & 10000                                                       & 10000                                                       \\ \hline
$N_{ic}$        & 10000                                                       & 10000                                                       & 10000                                                       & 30000                                                       \\ \hline
$\epsilon_{ic}$ & $5\times10^{-6}$                                            & $5\times10^{-6}$                                            & $5\times10^{-6}$                                            & $5\times10^{-6}$                                          \\ \hline
$xc_l$&$-3$&-3&-5&-9
          \\ \hline
$xc_r$&$3$&3&5&9
\\ \hline
\end{tabular}
    \caption{Parameters used by each method to solve the diffusion equation. $N_{int}$, $N_{bc}$ and $N_{ic}$ are the number of random points sampled at the interior, the boundary and at $t=0$, respectively. $xc_l$ and $xc_r$ are the limits of the computational domain used for each method.}
    \label{table:NN_parameters}
\end{table}

\begin{table}[hbt!]
    \centering
    \begin{tabular}{|c|c|}
    \hline
    Loss function threshold $\epsilon_k$     & Learning rate $\eta_k$ \\
    \hline
       - & 0.01 \\
       0.01  & 0.001\\
      0.005  & 0.0005\\
      0.001 & 0.0001\\
      0.0001   &0.00005 \\
    0.00001   & 0.00001\\
        \hline 
    \end{tabular}
    \caption{Learning rates as functions of the Loss function.}
    \label{table:learning_rates}
\end{table}

An important note when comparing these methods, is that they share only some parameters as their conditions for appropriate convergence may not be the same. In Table \ref{table:NN_parameters}, are shown some of the parameters used for the loss function minimization. In the table, the parameters $N_{int}$, $N_{bc}$ and $N_{ic}$ are the number of random points used at the interior of the domain, the points used at each boundary for the $x$ direction and the points used at the $x$ axis for $t=0$.

The parameter $\epsilon_{ic}$ plays an important role in the computation and follows the spirit of the subdomain splitting as will be discussed in the next section. In each of our formulations, when we minimize our solution trying to satisfy the initial condition, we do not use the exact initial condition. Instead, we minimize the loss function to approximate the exact initial condition until $A_L <\epsilon_{ic}$ is satisfied.  

Finally, as the loss function is reduced, the learning rate is also reduced. Initially, the learning rate is equal to $0.01$. However, when the value of the loss function reaches the value given in the left column (Table \ref{table:learning_rates}), it will take as the new learning rate the one given by the corresponding value on the right column. 

\subsection{Bistable Equation}
\begin{figure}[h]
	\centering
	\includegraphics[width=11cm, height=7.0cm]{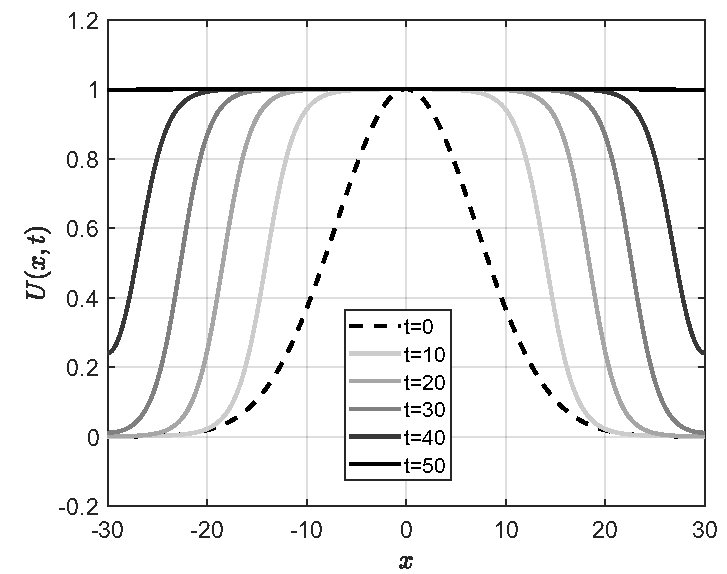}
	\caption{Numerical solution of the bistable equation obtained with the Chebyshev pseudospectral method given in \cite{oj22} with a fixed grid. }
	\label{Fig:bis_spec}
\end{figure}

This equation is given by 
\begin{equation}\label{Eq:bistable}
    u_t=Du_{xx}+u(u-a)(1-u).
\end{equation}
In order to test the methods we solve the equation with $D=1$ and $a=0.2$, for $x\in[-30,30]$, $t\in[0,50]$ and initial condition $u(x,0)=e^{-0.01x^2}$. This equation is solved for Neumann boundary conditions
 $\frac{\partial u}{\partial x}|_{x=x_L}=\frac{\partial u}{\partial x}|_{x=x_R}=0$. In this case, the analytical solution is not known but it can be obtained by some other methods like finite differences, pseudospectral methods or others. In Fig. \ref{Fig:bis_spec} we show the numerical solution obtained with the Chebyshev pseudospectral method given in \cite{oj22} with a fixed grid (parameters $N_i=900$, $N_c=3$ and $\Delta t=0.001$). This solution is computed at the fixed time values $t=0,5,10,15,20,25,30,35,40, 45$ and $50$. With the chosen number of subintervals $N_i$, collocation points per subinterval $N_c$, and time step $\Delta t$ the $l_{\infty}$ error is of the order of $10^{-7}$. This obtained solution is taken as the exact solution to make the comparisons of the developed methods in this work.

\begin{figure}[h]
	\centering
	\includegraphics[width=12cm, height=6.0cm]{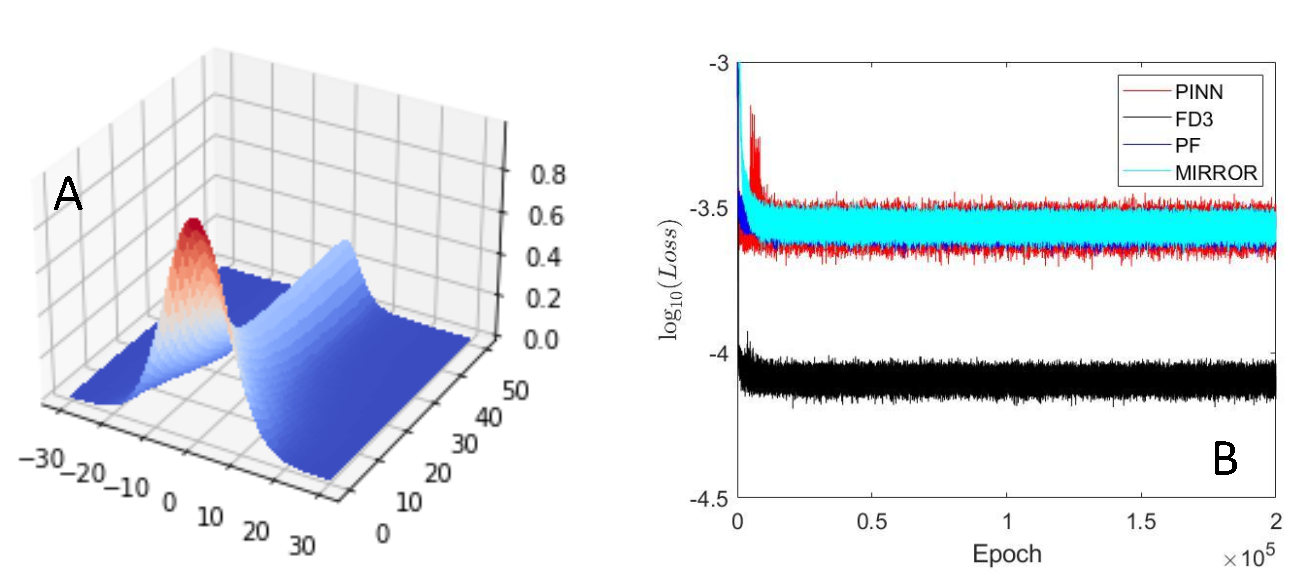}
	\caption{Numerical solutions of the bistable equation for the domain $x\in[-30,30]$ and $t\in[0,50]$ by each of the four methods. (A) The best approximation obtained with the PINN method after $2\times10^5$ epochs. The obtained solutions with the other three methods were practically the same. (B) The evolution of the loss versus the Epochs for each of the four methods.}
	\label{Fig:error_bist_nosubd}
\end{figure}

\begin{figure}[h]
	\centering
	\includegraphics[width=13cm, height=15.0cm]{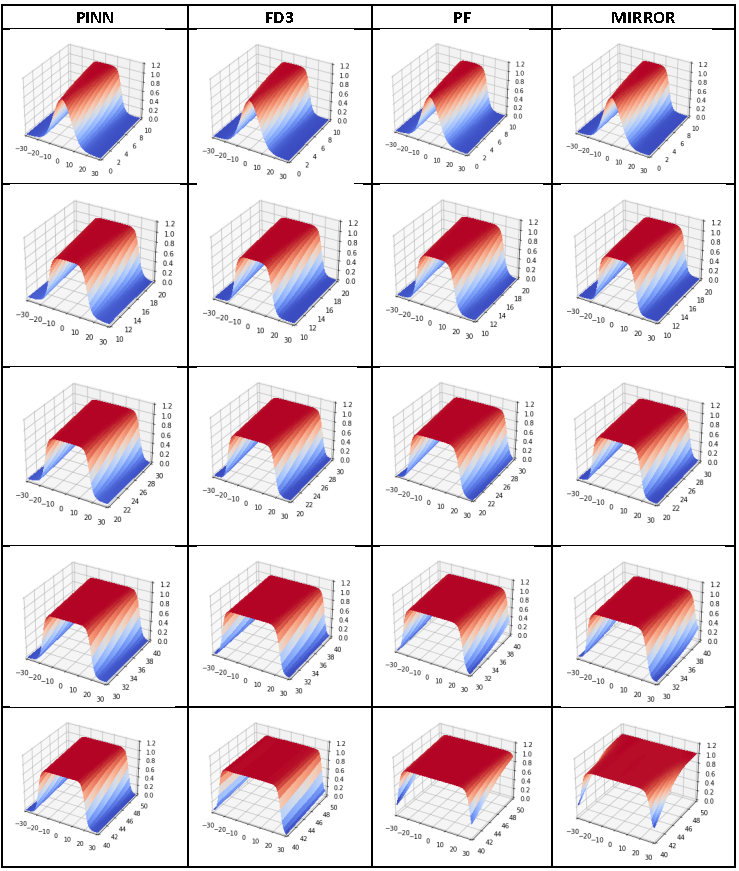}
	\caption{Numerical solutions of the bistable equation for the domain $x\in[-30,30]$ and $t\in[0,50]$. In the graph, the solutions are shown over a partition of $t\in[0,50]$ following the subdomain approach. The subdomains in time are $[0,10], [10,20], [20,30]$, $[30,40]$ and $[40, 50]$. }
	\label{Fig:Graficas_bistable}
\end{figure}

\begin{figure}[h]
	\centering
	\includegraphics[width=13cm, height=15.0cm]{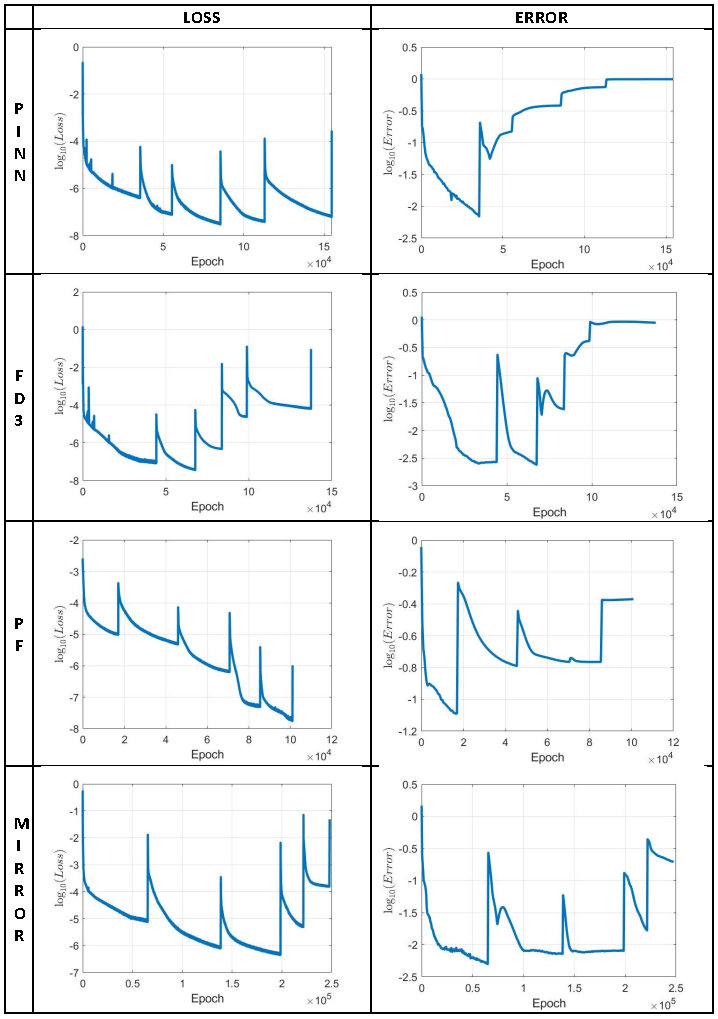}
	\caption{Loss and Error versus the Epoch, obtained with the four techniques. }
	\label{Fig:loss_error_bistable}
\end{figure}

The results of using the methods are shown in Figs. \ref{Fig:error_bist_nosubd}, \ref{Fig:Graficas_bistable} and \ref{Fig:loss_error_bistable}. In Fig. \ref{Fig:error_bist_nosubd} A it is shown the obtained solution with PINN after $2\times10(5)$ epochs. Clearly, the evolution of the propagating pulse does not follow the one obtained by the pseudospectral method in Fig. \ref{Fig:bis_spec}. In this scenario, it is observed that the pulse dies out. Similar solutions are obtained with the other three methods. On Fig. \ref{Fig:error_bist_nosubd}B it is shown the evolution of the loss function respect to the epochs for each of the four methods. Clearly, after $2\times10(5)$ epochs there is no improvement of the solution, which indicates that the method has ended up in a local minimum.  

The main issue is that the real solution has to steady states $U=0$ and $U=1$. In this scenario, a solution that remains close to $U=0$ or close to $U=1$ will provide regions in the parameter space that can lead to local minima. This observation was made in \cite{giampaolo2022} for the Gray-Scott equations. One possible way to provide with a better parameter search space is by considering simulations over smaller subdomains as described in Section \ref{Sec:subdomains}.
 
A trapping region was observed for this domain and larger ones in both, the $x$ and $t$ direction. For these larger domains, minimization of the loss function across the entire domain resulted in convergence to an incorrect local minimum that is difficult to escape. The implementation of subdomains helps to leave these minima and reach the correct solution.

For each of the methods, we have taken subintervals of length $10$ in the $t$ direction. Figure \ref{Fig:Graficas_bistable} illustrates the best solution obtained using the methods for the different subintervals. By eye comparison at $t=50$, from Figs. \ref{Fig:bis_spec} and \ref{Fig:Graficas_bistable}, we can see that the method's performance from worst to best is PINN, FD3, Phase Field and Mirror. Clearly, the solution at $t=50$ (last row in Fig. \ref{Fig:Graficas_bistable}), it should be very close to the constant value equal one, i.e.  $U(x,50)\approx 1$.
In this sense, for $t=50$, we observe that PINN and FD3 at the boundary remain close to zero, which is a bad approximation. Out of these two, the solution obtained with PINN will start rising at $t=60$ (Not shown), whereas for FD3 at $t=60$, will almost get to the value of $1$ (Not shown). The best solution is apparently obtained with the Mirror method.

To elucidate the performance of these methods, we show in Fig. \ref{Fig:loss_error_bistable}, left and right columns, the behavior of the loss and the error versus the epoch, respectively. From the right column, the order in which methods perform from worst to best is PINN, FD3, Phase Field, and Mirror. PINN performs well initially, but the solution worsens as the front reaches the boundary. FD3 performs better and again the solution is not good at the last two subdomains. Phase Field does not have the best accuracy at the beginning, but the error does not grow as for the first two methods. Finally, the Mirror technique is the one that shows better behavior in the evolution of the Error variable.

Observe that the behavior of the loss function for each method tells us the different nature of the algorithms to choose when they have to change the computations on the next domain. In the graph (Fig. \ref{Fig:loss_error_bistable}), the location at which the large amplitude peaks that are shown, represent the epoch at which the method changed the optimization procedure to the next subdomain. Observe that for each of the methods the value of the loss function is not a factor to decide when to state we have a good solution. For example, the Phase Field method reaches a value of almost $10^{-8}$ at the end of the computation, although, from the right column, the error is high. Clearly, there is no relation between the value of the loss function and the real error. Employing a fixed threshold on the loss function to determine transitions between optimization domains is demonstrably ineffective, as evidenced by the data in the figure (Fig. \ref{Fig:loss_error_bistable}). Observe that even within each of the four methods, the threshold value of the loss function for inducing a subdomain shift can vary.
 
The apparent advantage of the Mirror method could be questioned if we think about the time or epochs needed to obtain such a solution. From the left column of Fig. \ref{Fig:loss_error_bistable}, more iterations were needed to obtain the final solution for the Mirror method, compared to the other methods. Moreover, the parameters used for each of the methods are different. As the loss function to be minimized in each method is different, they may need different 
parameters in order to reach a good solution. Those parameters are shown in Table \ref{table:Bis_parameters}. The parameters for each method are such the methods were performing as their best with reasonable epochs. Observe that PINN, FD3 and PF does not coincide with parameters $N_{int}$, $N_{bc}$ and $N_{ic}$. However, particular runs were taken with equal values (Maximum value parameters) and the results were practically the same as those reported here. In the case of the Mirror technique, it was needed to add more nodes per layer due to the more complex solution to be approximated (Extended domain). Therefore, similar solutions were run with the first three methods in order to make the runs comparable. In all cases the Mirror method performed better. 

\begin{table}[]
\begin{tabular}{|l|l|l|l|l|}
\hline
                & PINN                & FD3                 & PF                  & MIRROR              \\ \hline
NN              & 4 Layers - 60 Nodes & 4 Layers - 60 Nodes & 4 Layers - 60 Nodes & 4 Layers - 90 Nodes \\ \hline
$N_{int}$       & 5000                & 10000               & 10000               & 10000               \\ \hline
$N_{bc}$        & 5000                & 5000                & 10000               & 10000               \\ \hline
$N_{ic}$        & 15000               & 15000               & 5000                & 9000                \\ \hline
$\epsilon_{ic}$ & $5 \times 10(-6)$           & $5 \times 10(-6)$            & $5 \times 10(-6)$            & $1 \times 10(-5)$            \\ \hline
$\epsilon_L$    & $1 \times 10(-7)$            & $1 \times 10(-7)$           & $1 \times 10(-5)$              &$1 \times 10(-6)$          \\ \hline
\end{tabular}
    \caption{Parameters used by the four methods to obtain the numerical solutions for the bistable equation (Eq. \ref{Eq:bistable}), using $D_L<\epsilon_L$ with $p=1000$ to change subdomains.}
    \label{table:Bis_parameters}
\end{table}

The experiment was repeated using $S_L$ given in Eq. \ref{Eq_Sl} as a changing domain criteria. By using the parameters given in Table \ref{table:Bis_parameters_SL}, we observed that the results basically, were the same as those shown in Table \ref{table:Bis_parameters}, but with larger amount of epochs. Observe that for these set of parameters, the value of $\epsilon_L$ for PINN and FD3 was taken smaller than the value for PF and Mirror. Even by reducing this value, we were not able to get a better solution than the ones obtained with PF and Mirror methods.

\begin{table}[]
\begin{tabular}{|l|l|l|l|l|}
\hline
                & PINN                & FD3                 & PF                  & MIRROR              \\ \hline
NN              & 4 Layers - 60 Nodes & 4 Layers - 60 Nodes & 4 Layers - 60 Nodes & 4 Layers - 90 Nodes \\ \hline
$N_{int}$       & 5000                & 5000               & 12000               & 15000               \\ \hline
$N_{bc}$        & 5000                & 5000                & 5000               & 5000               \\ \hline
$N_{ic}$        & 15000               & 15000               & 30000                & 45000                \\ \hline
$\epsilon_{ic}$ & $5 \times 10(-6)$           & $5 \times 10(-6)$            & $5 \times 10(-6)$            & $1 \times 10(-5)$            \\ \hline
$\epsilon_L$    & $5 \times 10(-8)$            & $5 \times 10(-8)$           & $5 \times 10(-6)$              &$5 \times 10(-7)$          \\ 
\hline
Epochs    & $3\times10(5)$            &     $3\times10(5)$       &  $2.7\times10(5)$             &  $4.1\times10(5)$     \\ 
\hline

\end{tabular}
    \caption{Parameters used by the four methods to obtain the numerical solutions for the bistable equation (Eq. \ref{Eq:bistable}), using $S_L<\epsilon_L$ criterion with $p=2000$ to change subdomains.}
    \label{table:Bis_parameters_SL}
\end{table}

\subsection{The Barkley's model}
\begin{figure}[h]
	\centering
	\includegraphics[width=11cm, height=7.0cm]{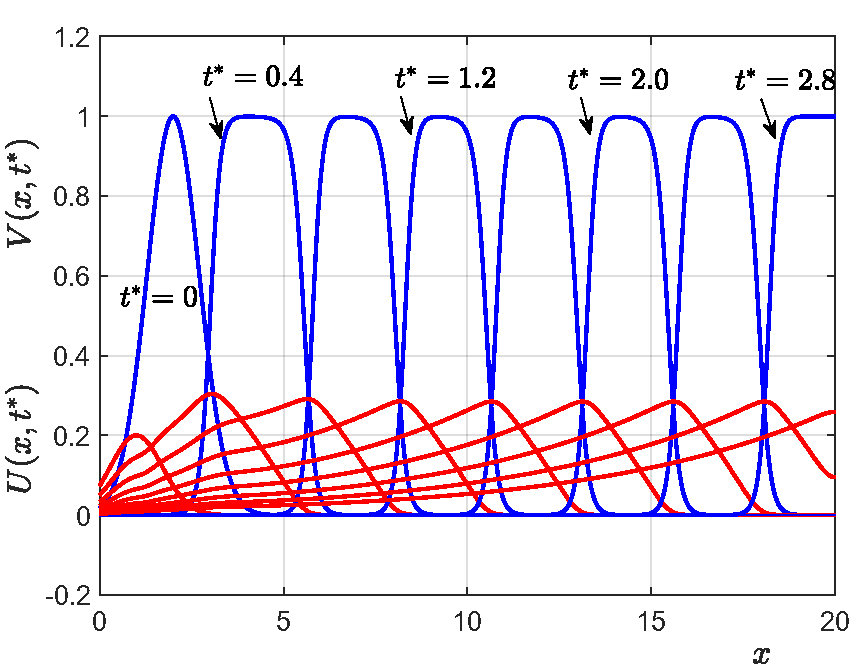}
	\caption{Numerical solution of the Barkley model obtained with the Chebyshev pseudospectral method given in \cite{oj22} with a fixed grid. $U(x,t^*)$ in blue, $V(x,t^*)$ in red.}
	\label{Fig:bark_spec}
\end{figure}

The last set of equations comes from the family of the Fitzhugh-Nagumo equations. It is known as the Barkley's model \cite{Barkley91} and is given by
\begin{equation}\label{barkeq}
\begin{array}{lcl}
u_t & = & \nabla^{2}u+\frac{1}{\epsilon}u(1-u)\left(u-\frac{(v+b)}{a}\right),  \\
v_t & = & u-v.
\end{array}
\end{equation}
In the equations, $u$ measures the membrane potential of an excitable cell and $v$ is a recovery variable that represents sodium channel reactivation and potassium channel deactivation \cite{Sherwood2014}. The solution is sought over the the spatial domain $x\in [0,20]$, integration time $t\in [0,3.2]$, initial conditions
\begin{equation}
\begin{array}{lll}
u(x,0)&=&e^{-(x-2)^2},\\
v(x,0)&=&0.2e^{-(x-1)^2},\\
\end{array}
\end{equation}
and no-flux boundary conditions for the $u$ variable. As the exact solution to this problem is not known, we build a good approximation ($U(x_i,t^*)\approx u(x_i,t^*)$, $V(x_i,t^*)\approx v(x_i,t^*)$) to the exact solution with the Chebyshev multidomain pseudospectral method given in \cite{oj22} with a fixed grid (parameters $N_i=2400$, $N_c=3$ and $\Delta t=1\times10^{-5}$). With these parameters, the obtained solution has an error of the order of $10^{-3}$. The solution is a pulse that propagates from left to right as shown in Fig. \ref{Fig:bark_spec}. The pulse arrives at the right boundary and then annihilates. The value of $U(x,3.2)\approx 0$. In the figure the solution is shown for different values of $t^*$.

For each of the four methods, we calculated the loss and the $l_{\infty}$ error as functions of the epochs to evaluate the performance. In order to calculate the $l_{\infty}$ error, we considered the exact solutions at times $t=0, 0.4, 0.8, 1.2, $ $1.6, 2, 2.4, 2.8$ and $3.2$. When minimizing in the $kth$ subinterval, the comparison to the exact solution was done with the solution obtained with the NN and the exact solution at a particular time which lied on that subinterval. The measurement of the error at these particular locations is just to provide with an idea of the performance of the method. For this particular example, we considered again Adam's method as an optimizer but with a fixed learning rate of $0.001$. 

\begin{figure}[h]
	\centering
	\includegraphics[width=11cm, height=15.0cm]{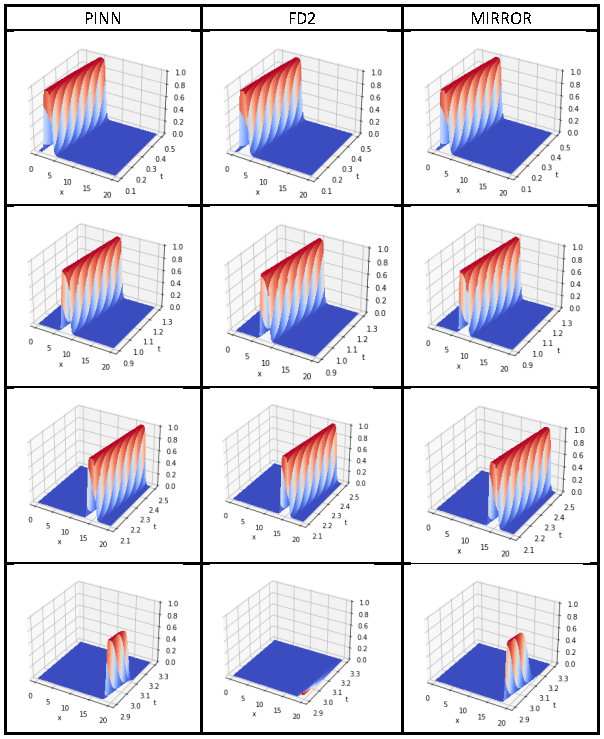}
	\caption{Numerical solution of the Barkley model obtained with the PINN, FD2 and Mirror methods.}
	\label{Fig:bark_sol_NN}
\end{figure}

In Fig. \ref{Fig:bark_sol_NN}, we show the solutions obtained with each of the developed methods. The domains at which minimization takes place are $[0,0.1]$, $[0.1,0.5]$, $[0.5,0.9]$, $[0.9,1.3]$, $[1.3,1.7]$, $[1.7,2.1]$, $[2.1,2.5]$, $[2.5,2.9]$ and $[2.9,3.3]$. Initially, a smaller time interval is taken of length $\Delta t=0.1$ as it helped to find the weights faster and not to lead to other minimum far from the solution. For each plot, it is shown the best approximation of the solution at that particular time interval. From the figure, it is clear that the PF method is the only one that do not performs well at all. The remaining three, provide with an appropriate solution at first sight. A better way to measure the performance is by looking at the evolution of the loss and the $l_{\infty}$ error as a function of the epochs.

\begin{figure}[h]
	\centering
	\includegraphics[width=8cm, height=15.0cm]{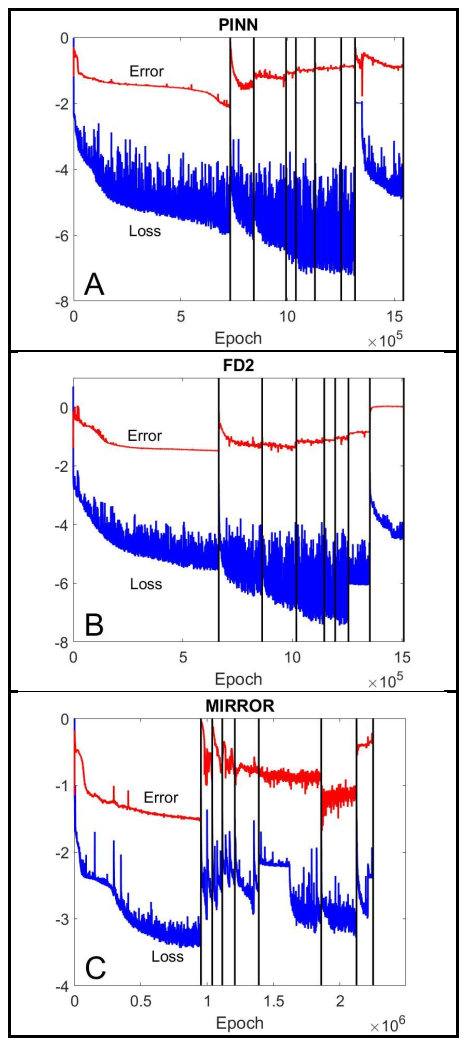}
	\caption{Evolution of the Error and the Loss for the PINN, FD2 and Mirror's methods.}
	\label{Fig:bark_errors_NN}
\end{figure}

In Fig. \ref{Fig:bark_errors_NN} we show the evolution of the loss and the $l_{\infty}$ error as a function of the epochs, for the PINN, FD2 and MIRROR's methods. The vertical lines in each of the figures, indicate the epoch at which the change of subomain takes place. Just at the left of the $kth$ vertical line the error plot provides the best reached error at the $kth$ subdomain. From the figure, we observe that both, the PINN and the Mirror method perform well, but PINN provides the answer with almost half the number of epochs. FD2, provides a reasonable answer until the last domain where the method loses the position of the pulse. An interesting phenomenon we observe for the Mirror method is that the error gets reduced as we move to the next subdomain. This is a very counter-intuitive result as we would expect that error gets accumulated as we move forward on domains. As a final comment is that PF method did not perform properly for this equation and therefore was not considered in the simulations.

\begin{table}[]
\begin{tabular}{|l|l|l|l|l|}
\hline
                & PINN                & FD2                 & PF                  & MIRROR              \\ \hline
NN              & 4 Layers - 60 Nodes & 4 Layers - 60 Nodes & 4 Layers - 60 Nodes & 4 Layers - 60 Nodes \\ \hline
$N_{int}$       & 500                & 500               & 12000               & 15000               \\ \hline
$N_{bc}$        & 500                & 500                & 5000               & 5000               \\ \hline
$N_{ic}$        & 500               & 500               & 30000                & 45000                \\ \hline
$\epsilon_{ic}$ & $1 \times 10(-6)$           & $1 \times 10(-6)$            & $5 \times 10(-6)$            & $1 \times 10(-5)$            \\ \hline
$\epsilon_L$    &  $5 \times 10(-9)$           & $5 \times 10(-9)$           &  $5 \times 10(-9)$             &$5 \times 10(-9)$          \\ 

\hline

\end{tabular}
    \caption{Parameters used by the four methods to obtain the numerical solutions for the Barkley's model (Eq. \ref{barkeq}), using $S_L<\epsilon_L$ criterion with $p=20,000$ to change subdomains.}
    \label{table:Barkley_parameters_SL}
\end{table}

\section{Conclusions and Discussion}\label{Sec:Conc_Disc}
In the present work, we have considered the problem of solving equations of the Reaction-Diffusion type in one dimension when no-flux boundary conditions were implemented. Four different techniques were presented and compared for three particular examples. A technique that splits the domain in the $t$ direction is also developed for long integration times. 

The domain split technique was mandatory to solve the equations for long integration times. Taking a subdomain with a large size in the time direction can result in the solutions converging to an incorrect local minimum, leading to an incorrect solution. However, dividing it into subdomains does not appear to be the final factor to obtain an appropriate approximation to the solution. Implementing boundary conditions differently compared to the traditional PINN implementation, can improve the approximation to the solution. Therefore, we can say that the two techniques, domain splitting and different ways of applying boundary conditions complement each other to obtain a desired solution. 

In this work, we observed that our techniques easily found the solutions of two of the examples considered, namely the bistable equations and the Barkley model. This is mostly due to the shape of the solutions which have a sigmoidal behavior.

For the presented methods, there is no clear way which technique performs better. It could be argued that the Mirror technique performs better, but it clearly has some disadvantages. It is more expensive than the rest of the methods as it requires a domain three times larger than PINN or the methods where the derivative condition is replaced by the finite difference approximation. Furthermore, when implementing the Mirror method in two dimensions, the equations can only be solved within rectangular domains. This is a very restrictive method. The other three methods can be easily extended to two dimensions but it may be possible that such methods become more expensive to obtain a proper solution. 

The Mirror technique is useful in the case of modeling propagating waves. For example, the speed and shape of a propagating pulse is not completely influenced by the boundary, if it is located far from it. However, the same is not true if we are looking to model a Turing pattern \cite{mu01}. In order to obtain a Turing pattern, the zero Neumann boundary condition plays an essential role and by imposing an artificial Dirichlet boundary condition over the extended domain, does not produce the desired patterns. 

As shown in this work, it is essential to implement a partition of the domain in the $t$ direction. However, the subdomain technique gives some challenges that need to be developed further. For example, depending on the equation to be solved, it may be possible to obtain a range of length values of the subdomains in the time direction that will provide with good results. For a given configuration of the neural network, and a particular RD equation, it should be possible to obtain a partition $t_0,t_1,t_2,...,t_N$ with $t_k=t_0+k\Delta t$, such that convergence to the solution occurs and at a relatively fast rate. Clearly, it should be easier to approximate the solution of a nonlinear RD equation over a smaller domain than over a large one. 

The subdomain technique can be implemented in more general ways. For example, the use of optimized weights in $\Omega_k$ as an initial guess to optimize over $\Omega_{k+1}$, can be changed to optimization over $\Omega_{k+1} \cup \tilde{\Omega}_k$ where $\tilde{\Omega}_k \subset\Omega_k$. 

For the finite differences approximations to the boundary conditions, it is clear that different approximations can be made. Second, third and fourth order were considered in this work. However, other approximations can be used as well, such as higher order finite differences or pseudospectral approximations.  

The methodology presented in this work has many limitations. For example, if we consider a single domain, once an initial condition has been approximated, does not follow automatically that we will get good approximations to the solution over the domain in the $(x,t)$ region. There can be transient behaviors in the solution that might not show up in the initial condition, which cannot be approximated with the same precision than the NN approximates the initial condition. This bad approximation, leads to get bad approximations to the solution. One possible way to avoid this problem, is to consider initial conditions with similar behavior than the propagating pulse profile, and by this avoiding any possible transient that could produce bad results. \\

In the worked examples we observed that the NN methodology in general solved the equations. However, the computing time for these methods is very high compared to the traditional methods of lines. However, it is worth to explore these new methodologies as not all the equations are easy to deal with, by traditional methods. On the other hand, there were equations that were not possible to solve appropriately. Such is the case of the Gray-Scott equations \cite{pet94}. In this scenario, the equations are highly nonlinear and pulse propagation follows as in the Barkley's model. However, for the Gray-Scott equations, pulses do not annihilate at the boundary but they get reflected. The solutions obtained with the methods in this work had a dispersive nature and the pulses die out. A possible solution involve subdomains with width in the $t$ direction very small, compared to the minimum integration time in which we can observe the physical phenomenon of reflection. This leads to a very expensive computation. Another option is to take large values of the number of nodes per layer. We used up to $140$ nodes per layer, and up to $7$ layers and still we were not able to obtain an appropriate answer. 

Finally, the choice of sampling points for each of the methods, as well as the used network were chosen by trial and error, such that the initial condition was approximated in a appropriate time and the solutions in the interior of the domain were obtained also with reasonable time.

\section{Acknowledgements}
The authors would like to acknowledge ACARUS at the University of Sonora, for their facilities for the numerical computations. JHJ was supported by NRF of Korea under the grant number
2021R1A2C3009648 and NRF2021R1A6A1A1004294412. JHJ was also supported partially by NRF grant by the Korea government (MSIT) (RS-2023-00219980).
%----------------------------------------------------------------

%\bibliographystyle{elsarticle-num}
%\bibliography{cas-refs}

\begin{thebibliography}{9}

\bibitem{brauer} F. Brauer, C. Castillo-Ch\'avez, Mathematical Models in Population Biology
and Epidemiology, Texts in Applied Mathematics, Springer, New
York, 2000.
\bibitem{epstein} I. R. Epstein, J. Pojman, Introduction to Nonlinear Chemical Dynamics,
Oxford University Press, Oxford, 1998.
\bibitem{keener_book} J. P. Keener, J. Sneyd, Mathematical Physiology, Interdisciplinary Applied
Mathematics, Springer, New York, 1998.
\bibitem{mu01} J. Murray, Mathematical Biology I and II, Interdisciplinary Applied
Mathematics, Springer, New York, 2002.
\bibitem{ty94} J. J. Tyson, What everyone should know about the Belousov-
Zhabotinsky reaction, in: S. Levin (Ed.), Frontiers in Mathematical
Biology. Lecture Notes in Biomathematics, Springer, Berlin, Heidelberg,
1994, pp. 569–587.
\bibitem{Anguelov2005} R. Anguelov, P. Kama, J.-S. Lubuma, On non-standard finite difference
models of reaction–diffusion equations, Journal of computational and
Applied Mathematics 175 (2005) 11–29. doi:10.1016/j.cam.2004.06.002.
\bibitem{Hu2012} G. Hu, Z. Qiao, T. Tang, Moving Finite Element Simulations for
Reaction-Diffusion Systems, Adv. Appl. Math. Mech. 4 (3) (2012) 365–
381.
\bibitem{os09} D. Olmos, B. D. Shizgal, Pseudospectral method of solution of the
Fitzhugh-Nagumo equation, Math and Comp. in Simul. 79 (7) (2009)
2258–2278.
\bibitem{oj22} D. Olmos-Liceaga, J.-H. Jung, A Chebyshev multidomain adaptive
mesh method for Reaction-Diffusion equations, Submitted (2021).
doi:10.48550/arXiv.2103.08484.
\bibitem{zhang2018} R. Zhang, Z. Wang, J. Liu, L. Liu, A compact finite difference
method for reaction–diffusion problems using compact integration factor
methods in high spatial dimensions, Advances in Difference Equations
2018 (274) (2018). doi:10.1186/s13662-018-1731-7.
\bibitem{Beck2021_2} C. Beck, S. Becker, P. Cheridito, A. Jentzen, A. Neufeld, Deep splitting
method for parabolic PDEs, SIAM J. Sci. Comput. 43 (5) (2021) A3135–
A3154. doi:10.1137/19M1297919.
\bibitem{Blechschmidt2021} J. Blechschmidt, O. G. Ernst, Three ways to solve partial differential
equations with neural networks — A review, GAMM - Mitteilungen
(2021). doi:10.1002/gamm.202100006.
\bibitem{Han2018} H. Han, A. Jentzer, W. E, Solving high-dimensional partial differential
equations using deep learning, PNAS 115 (34) (2018) 8505–8510.
doi:10.1073/pnas.1718942115.
\bibitem{Li2020} A. Li, R. Chen, A. Barati Farimani, Y. J. Zhang, Reaction diffusion
system prediction based on convolutional neural network, Scientific Reports
(10:3894) (2020) 1–9. doi:10.1038/s41598-020-60853-2.
\bibitem{ren2022} P. Ren, C. Rao, Y. Liu, J.-X. Wand, H. Sun, PhyCRNet: Physicsinformed
convolutional-recurrent network for solving spatiotemporal
PDEs, Comput. Methods Appl. Mech. Engrg. 389 (2022) 114399 (1–
21). doi:10.1016/j.cma.2021.114399.
\bibitem{scholz2022} C. Scholz, S. Scholz, Exploring complex pattern formation with convolutional
neural networks, American Journal of Physics 90 (2022) 141–151.
doi:10.1119/5.0065458.
\bibitem{Mishra2018} S. Mishra, A machine learning framework for data driven acceleration
of computations of differential equations, Mathematics in Engineering
1 (1) (2018) 118–146. doi:10.3934/Mine.2018.1.118.
\bibitem{raisi2018} M. Raissi, P. Perdikaris, G. Karniadakis, Physics-informed neural networks:
A deep learning framework for solving forward and inverse problems
involving nonlinear partial differential equations, Journal of Computational
Physics 378 (2019) 686–707. doi:10.1016/j.jcp.2018.10.045.
\bibitem{lagaris98} I. E. Lagaris, A. Likas, D. I. Fotiadis, Artificial Neural Networks for
Solving Ordinary and Partial Differential Equations, IEEE Trans. on
Neural Networks 9 (5) (1998) 987–1000. doi:10.1109/72.712178.
\bibitem{Cho2021} S. W. Cho, H. J. Hwang, S. Hwijae, Traveling wave solutions of partial
differential equations via neural networks, Journal of Scientific Computing
89 (21) (2021). doi:10.1007/s10915-021-01621-w.
\bibitem{Mattey2022} R. Mattey, S. Ghosh, A novel sequential method to train physics informed
neural networks for Allen Cahn and Cahn Hilliard equations,
Computer Methods in Applied Mechanics and Engineering 390 (2022)
114474. doi:10.1016/j.cma.2021.114474.
\bibitem{giampaolo2022} F. Giampaolo, M. De Rosa, P. Qi, S. Izzo, S. Cuomo, Physics-informed
neural networks approach for 1D and 2D Gray-Scott systems, Adv.
Model. and Simul. in Eng. Sci. 9 (5) (2022) 1–17. doi:10.1186/s40323-
022-00219-7.
\bibitem{Beck2021} C. Beck, M. Hutzenthaler, A. Jentzen, B. Kuckuck, An overview on deep
learning-based approximation methods for partial differential equations,
Cornell arxiv (2021). doi:10.48550/arXiv.2012.12348.
\bibitem{Meng2020} X. Meng, Z. Li, D. Zhang, G. E. Karniadakis, PPINN: Parareal
physics-informed neural network for time-dependent PDEs, Computer
Methods in Applied Mechanics and Engineering 370 (2020) 113250.
doi:10.1016/j.cma.2020.113250.
\bibitem{Zakeri2019} B. Zakeri, M. Khashehchi, S. Samsam, A. Tayebi, A. Rezaei, Solving
partial differential equations by a supervised learning technique, applied
for the reaction–diffusion equation, SN Applied Sciences 1:1589 (2019).
doi:10.1007/s42452-019-1630-x.
\bibitem{shekari2009} R. Shekari Beidokhti, A. Malek, Solving initial-boundary value problems
for systems of partial differential equations using neural networks and
optimization techniques, Journal of the Franklin Institute 346 (2009)
898–913.
\bibitem{Lu2021} L. Lu, X. Meng, Z. Mao, G. E. Karniadakis, DeepXDE: A Deep Learning
Library for Solving Differential Equations, SIAM Review 63 (1) (2021)
208–228. doi:10.1137/19M1274067.
\bibitem{Bueno2006} A. Bueno-Orovio, V. M. P´erez-Garc´ıa, F. H. Fenton, Spectral methods
for partial differential equations in irregular domains: The spectral
smoothed boundary method, SIAM J. Sci. Comput. 28 (3) (2006) 886–
900. doi:10.1137/040607575.
\bibitem{Barkley91} D. Barkley, A model for fast computer simulation of waves in excitable
media, Phys. D: Nonlinear Phenomena 49 (1-2) (1991) 61–70.
doi:10.1016/0167-2789(91)90194-E.
\bibitem{Sherwood2014} Sherwood. W. E., Fitzhugh–nagumo model, in: D. Jaeger, R. Jung (Eds.), Encyclopedia
of Computational Neuroscience, Springer, New York, NY,
2014, pp. 1–11.
\bibitem{pet94} V. Petrov, S. Scott, K. Showalter, Excitability, wave reflection and wave
splitting in a cubic autocatalysis reaction diffusion system, Phil. Trans.
R. Soc. Lond. A 347 (1994) 631–642.

\end{thebibliography}

%-------------------------------------------------------------------------------

\end{document}